\documentclass[manmat,referee]{svjour}
\usepackage{amsmath,amssymb}
\author{ Rosa Cid-Mu$\tilde n$oz\and  Manuel
Pedreira}
\title{Classification of Incidence
Scrolls(I)}
\date{}
\institute{\emph{Authors' address:} Departamento de
Algebra, Universidad de Santiago de Compostela, $15782$
Santiago de Compostela, Galicia, Spain. Phone:
34-81563100-ext.13152. Fax: 34-81597054.{\tt e}-mail: {\tt
rosacid@usc.es}; {\tt
pedreira@zmat.usc.es}}
\newtheorem{teo}[subsection]{Theorem}
\newtheorem{defin}[subsection]{Definition}
\newtheorem{prop}[subsection]{Proposition}
\newtheorem{cor}[subsection]{Corollary} 
\newtheorem{lemma1}[subsection]{Lemma}
\newtheorem{example1}[subsection]{Example}
\newtheorem{rem}[subsection]{Remark}

\newcommand\om{\Omega}
\newcommand\p{{\bf P}}
\newcommand\E{{\cal E}}
\newcommand\Te{{\cal O}}

\newcommand\X{{\cal X}}
\newcommand\lu{{\cal L}}
\font\euf=eufm10 at 14pt
\newcommand\e{\mbox{\euf e}}
\def\d{\mbox{\euf b}}
\def\de{\mbox{\euf d}}
\newcommand\be{{\cal B}}

\def\qed{\hspace{\fill}$\rule{1.5mm}{1.5mm}$}
\newcommand\lrw{\longrightarrow}
\newcommand\rw{\rightarrow}
\newcommand\sub{\subset}

\begin{document}
\combirunning{Rosa Cid, M. Pedreira}
\maketitle
\begin{abstract}
The aim of this
paper is to obtain a classification of scrolls of genus $0$
and $1$, which are defined by a one-dimensional family of
lines meeting a certain set of linear spaces in $\p ^n$.
These ruled surfaces  will be called incidence scrolls and
such a set will be the base of the incidence
scroll. Unless otherwise stated, we assume that the
base spaces are in general position.
\end{abstract}

{\Large\bf Introduction}

Throughout this paper, the
base field for algebraic varieties is
${\mathbb{C}}$. Let $\p^n$ be the n-dimensional
complex projective space and $G(l,n)$ the
Grassmannian of $l$-planes in
$\p^n$. Then $R^d_g \sub \p^n$ denotes a scroll of degree $d$
and genus $g$. We will follow the notation of \cite{hatshor}.

Projective ruled surfaces with the
property that the generators meet a suitable number of
linear spaces in
$\p^n$ are known classically. The first examples of this situation 
are the smooth quadric surface in $\p^3$ and the elliptic quintic
scroll in
$\p^4$ formed by the lines which meet five planes
in general position. Although this is not so in all the ruled
surfaces, there does not exist a classification of such
ruled surfaces. Our purpose is to investigate the scrolls in
$\p^n$ whose generators meet a certain set
of linear spaces in general position. The ruled surfaces so 
obtained will be called incidence scrolls, and such a
set will be a base of the incidence scroll. 

It is useful to represent a scroll in $\p^n$ by a curve
$C\sub G(1,n)\sub \p^N$. The lines which intersect
a given projective space $\p^r$ in $\p^n$ are represented by
the points of the special Schubert variety $\om
(\p^r,\p^n)$. Each
$\om (\p^r,\p^n)$ is the intersection of $G(1,n)$ with a
certain linear space in $\p^N$. Since $G(1,n)$ has dimension
$2n-2$ and we need a curve, we must impose $2n-3$ linear  
conditions on $G(1,n)$. Consequently, the choice of linear
spaces is not arbitrary. Any set of linear spaces in $\p^n$
which imposes $2n-3$ linear conditions on $G(1,n)$ is the base
of a certain incidence scroll. The background about Schubert
varieties can be found in \cite{Kleiman}.

We will expose a method to know when a scroll is determined
by incidences. The affirmative solution would allow us to
obtain a base for each incidence scroll (there is really
only one way to choose this base). This is
possible because the families of directrix curves provide a
natural and intrinsic characterization of the incidence
scrolls.

Our first step will be to give some general properties of
ruled surfaces. For more details we refer the reader to
\cite{hatshor}. We do it in section \ref{1}, where we define a
particular class of surfaces, the scrolls. 
In section \ref{2}, having introduced the notion of incidence
scroll, we have compiled some basic properties of such a
scroll. The degree of the scroll given by a general base is
provided by Giambelli's formula. Our main contribution is
to establish $(IS)$ which together with \ref{25} - \ref{l62}
gives the first base spaces of the incidence scroll. 

Section \ref{3} is devoted to the study of deformations of a given
incidence scroll. These form a powerful tool for obtain a
finite family of incidence scrolls from a given one. If the
incidence scroll $R^d_g \sub \p^n$ breaks up into
$R^{d_{1}}_{g_{1}} \sub \p^r$ and $R^{d_{2}}_{g_{2}} \sub
\p^s$ with
$\delta$ generators in common, then $d=d_{1}+d_{2}$ and $
g=g_{1}+g_{2}+\delta-1$. 

Sections \ref{4} and \ref{seiss} contain our main results about
rational and elliptic incidence scrolls. In the former
section, we show that the only rational scrolls determined by
incidences are those of general type and those with a
directrix line. In the latter section, we give all the
elliptic incidence scrolls. In each case, we find the base
of the scroll. Accordingly, we will apply the following
results to obtain a large number of incidence scrolls: 
\begin{itemize}
\item Bertini's theorem: if $\wedge$ is a
reduced linear system on a complete nonsingular variety $V$
such that $dim(\Phi_{}(V))\geq 2$, then the generic element
of $\wedge$ is irreducible (see \cite{itaka}, p. 252).
\item A family of closed subschemes $f\colon X\sub T\times
\p^n \lrw T$ is flat if and only if the Hilbert polynomials of
the fibres are the same and hence also the arithmetic genus.
Every algebraic family of normal varieties parametrized by
a nonsingular curve $T$ over an algebraically closed field
is a flat family of schemes. For any flat family, having a point
$t\in T$ with $X_t$ birationally equivalent to $X_0$, we say
that
$X$ is a (global) deformation of $X_0$ (see \cite{hatshor},
pp. 256-266).
\end{itemize}

We then prove the main theorem of this sections, which
determines all the incidence scrolls of genus $0$ and $1$.
\newline
{\bf Theorem}(see Theorems \ref{cla} and \ref{c.e})
{\it Let $X$ be a ruled surface over the curve
$C$ of genus $g=0,1$. Let
$H\sim C_o+\d f$ be a very ample divisor on
$X$ with $m=deg(\d)$ and let {\small $\Phi_H:X \rw
R_g^{2m-e}
\sub
\p^{2m-e-2g+1}$} be the closed immersion defined by
$|H|$. Then $R_g^{2m-e}$ is an incidence scroll
if and only if it satisfies one of the following
conditions:
\begin{enumerate}
\item $g=0$ and $e=0,1$;
\item $g=0$ and $m=e+1$;
\item $g=1, \, \,e=-1$ and $m=2$;
\item $g=1, \, \,\e\sim 0$ and $m=4$;
\item $X$ decomposable, $g=1, \, \,0 \leq e\leq 3$
and $m=e+3$.
\end{enumerate}
 }

It is to be expected that the treatment of general case will
be answered shortly. A paper on incidence scrolls in
$\p^n$ is to be written. The results on this paper 
belong to the Ph.D. thesis of the first author whose advisor 
is the second one.

We thank the referee for many detailed remarks of
a first version of this paper.

\section{Ruled Surfaces} \label{1}

A ruled surface is a surface $X$ together with a
surjective morphism $\pi \colon X\lrw C$ to a (connected 
nonsingular) curve $C$ such that the fibre $X_y$ is
isomorphic to $\p^1$ for every point $y\in C$, and such
that $\pi$ admits a section. There exists a locally free
sheaf $\E$ of rank 2 on $C$ such that 
$X\cong \p (\E)$ over $C$. Conversely, every such $\p
(\E)$ is a ruled surface over $C$. Let $\Te_X (1)$ be the 
invertible sheaf $\Te _{\p(\E)}(1)$. Then there is an one-to-one
correspondence between sections $\sigma \colon C\lrw X$ and surjections 
$\E \lrw \lu \lrw 0$, where $\lu$ is an invertible sheaf on $C$, given by
$\sigma^{\star}(\Te_X(1))$. Moreover, if $D$ is any section
of $X$, corresponding to a surjection $\E \lrw \lu \lrw 0$,
and if $\lu = \lu (\de )$ for  some divisor $\de$ on $C$, then
$deg(\de)=C_o \cdot D$, and $D$ and $C_o+(\de -\e )f$ are linearly 
equivalent, written $D\sim C_o+(\de -\e )f$.

It is possible to write $X \cong \p (\E)$ where $\E$ is a
locally free sheaf on $C$ with the property that
$H^0(\E)\not= 0$ but for all invertible sheaves $\lu $ on $C$
with $deg\lu <0$, we have $H^0(\E \otimes \lu )=0$. In this
case we say $\E$ is normalized. Let $\e$  be the divisor on $C$
corresponding to the invertible $\bigwedge ^2\E$ and let $e =-deg\e$.
Furthermore there is a section $\sigma _0 \colon C \lrw X$
with image $C_o$ such that $\lu (C_o)\cong \Te_X(1)$. Fix
such a section $C_o$ of $X$. If $\d$ is any divisor on $C$,
then we denote the divisor $\pi^\star \d $ on $X$ by $\d f$,
by abuse of notation. Thus any element of Pic($X$)
can be written $aC_o+\d f$ with $a\in {\mathbb{Z}}$ and
$\d
\in $ Pic($C$).

Let $X$ be a ruled surface over the curve $C$ of genus $g$, determined by a
normalized locally free sheaf $\E$. If
$\E$ is decomposable, then
$\E \cong \Te_C\oplus \lu$ for some $\lu$ with $deg\lu \leq 0$. All
values of $e\geq 0$ are possible. If $\E$ is indecomposable,
then $-g\leq e\leq 2g-2$.

A scroll is a ruled surface embedded in $\p^N$
in such a way that the fibres $f$ have degree 1. If we take a
very ample divisor on $X$, $D\sim aC_o+\d f$, then the
embedding $\Phi_D \colon X \lrw R^d_g\sub
\p^N=\p(H^0(\Te_X(D)))$ determines a scroll when $D\cdot
f=a=1$. Throughout this work, we will restrict our attention
to very ample divisors so that $\Phi_D(X)$ be a scroll.

\section{Definition of Incidence Scroll} \label{2}

\begin{defin} \label{de1}
{\em A scroll $R\sub \p^n$ is said to be {\it an incidence
scroll} if $R$ is generated by the lines which meet a certain
set $\be$ of linear spaces in $\p^n$. Such a set is called a
{\it base} of the incidence scroll.}
\end{defin}

Since a scroll $R^d_g\sub \p^n$ is represented by a curve
$C^d_g\sub G(1,n)$, there is another definition
of incidence scroll which is equivalent to the previous
one.

\begin{defin} \label{de2}
{\em A scroll $R\sub \p^n$ is {\it an incidence scroll}
if the correspondent curve in $G(1,n)$ is an intersection of
special Schubert varieties $\om(\p^{r}, \p^{n})$, $0\leq r
<n-1$.}
\end{defin}  

Such a base will be denoted by:
$$\be =\{\p^{n_1},\p^{n_2},\cdots ,\p^{n_r}\}$$ 
or we will write it simply $\be$ when no confusion can arise,
where $n_1\leq n_2\leq \cdots \leq n_r$, i.e.,
$C=\bigcap_{i=1}^{r}(\om(\p^{n_i}, \p^{n}))\sub G(1,n)$. If
$n_1=\cdots=n_i, \, \, n_{i+1}=\cdots=n_j,\, \,\dots \, \,
\mbox{and}\,\, n_{j+1+k}=\cdots=n_r$, then we write it $\be
=\{i\,\p^{n_1},(j-i)\p^{n_j},\cdots , (r-j-k)\p^{n_r}\}$,
for short.

We can certainly assume that $n\geq 3$. This involves no
loss of generality because the only scroll in $\p^2$ is
$\p^2$. The plane can be obtained as incidence scroll in
$\p^n$ if we take
$\be =\{\p^{0},(n-2)\,
\p^{n-2}\}$. If we want to see it as nondegenerate surface,
we must work in $\p^2$ with $\be=\{\p^{0}\}$ (i.e.,
$\p^1=\om(\p^{0},\p^2)$).

We are  interested in incidence scrolls where the base is
formed by linear spaces in general position. By
general position we will mean that
$(\p^{n_1},\cdots ,\p^{n_r})\in {\cal W}=G(n_1,
n)\times \cdots \times G(n_r, n)$ is contained in
a nonempty open subset of
${\cal W}$. Therefore, unless otherwise stated, we will
work with general base spaces. For simplicity of
notation, we abbreviate it to base in general
position. 

\begin{teo} \label{S} The intersection $\bigcap_{j=1}^m
\Omega(\p^{n_j},\p^{n})\sub G(l,n)$ of the special Schubert
varieties as\-soci\-ated to linear subspaces
$\p^{n_j},\, j=1,\cdots,m$, of dimension
$n_j$ such that $(l+1)(n-l)-\sum_{j=1}^{m}(n-n_j-1) >0$, is
connected. Moreover, the intersection is irreducible of
dimension
$(l+1)(n-l)-\sum_{j=1}^{m}(n-n_j-1)$ for a general choice
of the subspaces $\p^{n_j}$.
\end{teo}
{\bf Proof.} See \cite{sols}, Theorem 1.1.

With the previous theorem, we can prove immediately an
important consequence which provides a characterization of
incidence scrolls.

\begin{prop} \label{IS} The intersection
$C=\bigcap_{i=1}^{r}
\Omega(\p^{n_i},\p^{n})$ of special Schubert varieties
associated to linear spaces $\p^{n_i},\, i=1,\cdots,r$ in
general position is a conex irreducible curve of
$G(1,n)$ if and only if it verifies the following
equality
$$ rn-(n_1+n_2+\cdots +n_r)-r=2n-3 \eqno(IS) $$
\end{prop} 
{\bf Proof.} $''\Rightarrow''$ Since $C\subset
G(1,n)$ is a conex irreducible curve, we have
$1=2(n-1)-rn+r+\sum_{i=1}^{r} n_i\, \, \mbox{(see Theorem
\ref{S})}
\iff rn-\sum_{i=1}^{r} n_i -r=2n-3$.

$''\Leftarrow''$ Since $2(n-1)-\sum_{i=1}^{r} n-n_i-1 >0$,
we conclude that
$C$ is connected, by Theorem \ref{S}. Moreover, we easily
see from Theorem \ref{S} that it is irreducible of
dimension $1$ for a general choice of the $\p^{n_i},\,
i=1,\cdots,r$. \qed

Each $\om(\p^{n_i}, \p^n)$ is the intersection of
$G(1,n)$ with a linear space in
$\p^N$. Moreover, such a linear space is a hyperplane if
and only if we have $n_i=n-2$. If every
dimension of the base spaces is equal to
$n-2$, then we will need exactly $(2n-3)\,
\p^{n-2}$'s for obtain a curve in $G(1,n)$.

Some relevant properties of incidence scrolls will be
indicated. These are elementary but of great importance if
we want to develop a theory about incidence scrolls.

\begin{enumerate}
\item  A hyperplane does not impose
conditions on $G(1,n)$ because
$\om(\p^{n-1},\p^{n})$ $=G(1,n)$. Then we can
assume $n_r\leq n-2$. 
\item If $n_i+n_j<n-1$ then the scroll is degenerated. Under
the above assumptions, $\om(\p^{n_i}, \p^n)
\cap \om(\p^{n_j}, \p^n)=\om(\p^{n_i}, \p^{n_i+n_j+1})
\cap \om(\p^{n_j}, \p^{n_i+n_j+1})$.
\item A nondegenerate (irreducible) incidence scroll
cannot have double points. In other case, if $P$
is a double point, there are two generators of the
scroll $g_1, \,g_2$ through $P$. Thus, the
incidence scroll is not irreducible, containing the plane
pencil determined by them. In $G(1,n)$, we have
$\bigcap_{i=1}^{r}
\om(\p^{n_i},
\p^n)=\om(P, \p^2) \bigcup(\bigcap_{j=1}^{s} \om(\p^{m_j},
\p^m))$.  
\item The generators of the scroll meet each base space in
points of a directrix curve.
\item The degree of the scroll can be
calculated by a one-to-one correspondence between any two of the
directrix curves. It can also be computed by
Giambelli's formula, i.e., it is the number of lines in
$\p^n$ which intersect $\p^{n_1},\p^{n_2},\cdots ,\p^{n_r}$
and a generic $\p^{n-2}$. If $deg(R)=d$, then we have the
following equality of Schubert cycles: $\om(n_1, n) \cdots
\om(n_r, n)=d\om(0, 2)$. 
\end{enumerate}

From now on, we will talk about a decomposable
incidence scroll if the corresponding ruled surface 
$X=\p(\E)$ has $\E$ decomposable. If $\E$ is indecomposable, 
then we will talk of an indecomposable incidence scroll.
 
Let $X=\p(\E)$ be a ruled surface
over the curve $C$ of genus $g$, determined by a
decomposable normalized bundle $\E\cong
\Te_C\oplus\Te_C(\e)$ such that $deg(-\e)=e\geq
0$. Let $H\sim C_o+\d f$ be the very ample divisor
on $X$ with $m=deg(\d)$ which gives the immersion
of the ruled surface as the scroll $R^d_g \sub
\p^n$ such that $d=2m-e$ and $n=2(m-g)-e+1+i$, being $i$ the 
speciality of the
scroll ($i=h^1(\Te_C(\d))+h^1(\Te_C(\d +\e))$).
Geometrically, $X$ has two disjoint directrix,
denoted by $C_o$ and $C_1$, such that $C_1\sim
C_o-\e f$. Moreover, it is easy to check that
these satisfy $\phi_{\d + \e }
\colon C_o\lrw C_g^{m-e} \sub \p^{m-e-g+i_1}$ and 
$\phi_{\d } \colon C_1\lrw C_g^{m}\sub
\p^{m-g+i_2}$, being $i_1=h^1(\Te_C(\d +\e))$ and
$i_2=h^1(\Te_C(\d))$. Therefore
$i=i_1+i_2$ and $\p^{m-e-g+i_1}\cap \p^{m-g+i_2} =
\emptyset$.

\begin{prop} \label{25} If $R^d_g \sub \p^n$ is a
decomposable incidence scroll with base in general
position, then
$\p^{m-e-g+i_1}$ and $\p^{m-g+i_2}$ are base spaces.
\end{prop}
{\bf Proof.} Suppose $\be$ a base with
$\sum_{i=1}^{r}(n-n_i-1)=2n-3$ such that
$\p^{n_i}$'s, $i=1\cdots,r$, are in general
position and
$m-e-g+i_1< n_i \leq n-2$. Since the scroll has a minimum
directrix curve $C^{m-e}_g \sub \p^{m-e-g+i_1}$, we can
take a generic hyperplane $\p^{n-1}$ through
$\p^{m-e-g+i_1}$. Then there are $m$ lines in
$\p^{n-1}$ which meet $\p^{{n_1}-1}, \p^{{n_2}-1},
\cdots $ and $\p^{{n_r}-1}$. Since
each $\p^{{n_i}-1}$ imposes
$n-n_i-1$ independent conditions on $\p^{n-1}$,
it is impossible.

$\p^{m-e-g+i_1}$ imposes $m-g+i_2$ conditions on
$G(1,n)$. Then we can  consider more subspaces to form
the base. There is a $\p^{m-g+i_2} \in \be$. For
otherwise, suppose $\be
=\{\p^{m-e-g+i_1},\p^{n_1}, \cdots ,\p^{n_r}\}$
(in general position) with
$m-g+i_2<n_i\leq n-2$. Then we can take a generic
hyperplane
$\p^{n-1}$ such that $\p^{m-g+i_2} \sub \p^{n-1}$
and arguments similar to the above imply that
$\be$ is not the base of
$R_g^d\sub \p^{n}$.\qed

\begin{teo} \label{26}
Let $R_g^d\sub \p^{n}$ be an incidence scroll with base
in general position. Then:$$ R\quad is \quad 
decomposable
\iff
\p^{n_1}\cap \p^{n_2}=\emptyset.$$
\end{teo}
{\bf Proof.}''$\Leftarrow$'' Since $\p^{n_1}\cap
\p^{n_2}=\emptyset$, the directrix curves $C^{d_1}_g\sub
\p^{n_1}$ and $C^{d_2}_g\sub\p^{n_2}$ satisfy $C^{d_1}_g\cap
C^{d_2}_g=\emptyset$, i.e., $R^d_g$ is decomposable.

''$\Rightarrow$'' By Proposition \ref{25}, there are
two directrix curves $C^{d_1}_g\sub\p^{n_1}$ and
$C^{d_2}_g\sub\p^{n_2}$ such that $C^{d_1}_g\cap
C^{d_2}_g=\emptyset$. Suppose
$n_j=d_j-g+i_j$ for $j=1,2$. Since
$n=n_1+n_2+1$ and
$\p^{n_1}$ and $\p^{n_2}$ are in general position,
$\p^{n_1}\cap
\p^{n_2}=\emptyset$.\qed

If $R$ is decomposable, $\p^{m-e-g+i_1}$ and
$\p^{m-g+i_2}$ impose $2(m-g)-e+i$ conditions on
$G(1,n)$. Since
$2n-3-(2(m-g)-e+i)=n-2\geq 1$, we must consider
more subspaces to form the base.

In particular, if $e \leq 2g-2$, then
$h^0(\Te_X(C_o-\e f))=e+2-g$. Moreover, if $g=0,1$,
then we can show (by Bertini's theorem) that there
are
$e+2-g$ linearly independent directrix curves
$C\sim C_o-\e f$, which generates $|C_o-\e f|$.

\begin{prop} \label{l62} Let $R^{2m-e}_g\sub
\p^{2(m-g)-e+1+i}$ be a decomposable incidence
scroll with base in general position. 
\begin{enumerate}
\item[(a)] If $\e \not\sim 0$, then there are 
$(e+2-g)\, \p^{m-g+i_2}$'s in
$\be$, whenever possible, i.e., when
$m-g+i_2+(e+2-g)(m-e-g+i_1)\leq 4(m-g)-2e+2i-1$
and the generic curve in $|C_o-\e f|$ is
irreducible.
\item[(b)] If $\e \sim 0$, then there are 
$3\, \p^{m-e-g+i_1}$'s in
$\be$.
\end{enumerate}
\end{prop}
{\bf  Proof.} $(a)$ We proceed by induction in
$j=e+2-g$. The proposition is true for $j=1$, by
Proposition \ref{25}. Supposing the proposition 
true for
$j$, we prove it for $(j+1)\,
\p^{m-g+i_2}$'s. Let $\be =\{\p^{m-e-g+i_1},
\p^{n_1},$ $\p^{n_2},\cdots ,\p^{n_r}\}$ be a base of the
scroll in general position with
$m-g+i_2=n_1= \cdots =n_j<n_{j+1}\leq \cdots
\leq n_r\leq 2(m-g)-e+i-1$. Since the scroll
has at least one directrix curve $C^m_g \sub
\p^{m-g+i_2}\notin
\be$, which is linearly independent from the other
directrix curves contained in base spaces of
dimension $m-g+i_2$, we can take a generic
hyperplane
$\p^{2(m-g)-e+i}$ through $\p^{m-g+i_2}$. Then
there are
$m-e$ lines in $\p^{2(m-g)-e+i-j}$ which meet
$\p^{m-e-g+i_1-1},
\p^{n_1-j},$
$\cdots,\p^{{n_j}-j} ,
\p^{n_{j+1}-j-1},$ $\cdots $ and $\p^{{n_r}-j-1}$,
which is impossible.

The proof for $(b)$ is similar.
\qed

\section{Degeneration} \label{3}

\begin{prop} \label{123}
Let $R^d_g \sub \p^n$ be an incidence scroll with base 
$\be$ in general position. Suppose that $\p^{n_i}$ and
$\p^{n_j}$ lie in a  hyperplane $\p^{n-1}$ and have in
common
$\p^{m}$, 
$m=n_i+n_j-n+1$. Then
the scroll breaks up into: 
\begin{enumerate} 
\item[-] $R^{d_1}_{g{_1}} \sub \p^{n}$
with base $\dot \be
=\{\p^{m},\p^{n_1},\cdots,\widehat{\p^{n_i}},\cdots,\widehat{\p^{n_j}},
\cdots,\p^{n_r}\}$ (which is possibly degenerate);
\item[-] $R^{d_2}_{g_2} \sub \p^{n-1}$
with base $\ddot \be =\{\cdots ,
\p^{n_{i-1}-1},\p^{n_i},\p^{n_{i+1}-1},\cdots
,\p^{n_{j-1}-1},\p^{n_j},$ $\p^{n_{j+1}-1},\cdots \}$  
\end{enumerate}
which have $\kappa \geq 1$ generators in common. Then,
$d=d_1+d_2$ and $g=g_1+g_2+ \kappa -1$.

Moreover, if $m=0$, then the incidence scroll breaks up
into a plane and an incidence scroll $R^{d-1}_g \sub
\p^{n-1}$ with base $\ddot \be$ in general position. 
\end{prop}
{\bf Proof.} By assumption we have a scroll formed by the
lines which pass through $\p^m$ and a scroll generated by
the lines which intersect the base spaces and lie in
$\p^{n-1}$. The three scrolls are represented by curves
$C=\bigcap_{k=1}^{r} \om(\p^{n_k}, \p^n)$, 
$\dot C=\om(\p^{m}, \p^n)\cap (\bigcap_{k\not=i,j}
\om(\p^{n_k}, \p^n))$ and $\ddot C=\om(\p^{n_i},
\p^{n-1})\cap
\om(\p^{n_j},\p^{n-1})\cap(\bigcap_{k\not=i,j}
\om(\p^{n_k-1}, \p^{n-1} ))$ because 
every one of them satisfies $(IS)$.

Write $A_i=\om(\p^{n_i}, \p^n)), \, i=1,\cdots,r$, where
the $\p^{n_i}$'s are in general position and
$B_j=\om(\p^{n_j}, \p^n)), \, j=1,\cdots,r$, where
$\p^{n_1}\cap \p^{n_2}=\p^m$ and the other subspaces are in
general position. Then there is an invertible linear
transformation of $\p^N$ into itself which carries
$G(1,n)$ into itself and $A_i$ into $B_i$ (see
\cite{Kleiman} , Proposition 4). Consequently, we
conclude that $A_i$ and $B_i, \, \, i=1,\cdots
,r$, define the same cohomology class. In this way
we obtain what the set of lines which
simultaneously intersect the subspaces of
$\be$ can be continuously deformed into the union of two
set of lines. A set formed by the lines which pass through
a fixed point and lie in a plane and other set formed
by the lines which simultaneously intersect the subspaces of
$\ddot \be$. Since the various subvarieties in a continuous
system are all assigned the same cohomology class, we find
the following equality of Schubert cycles, $C=\dot C+\ddot
C=d\om(0, 2)$. So,
$d=d_1+d_2$.

We have in $G(1,n)$ a short exact sequence: $$
0\rw \Te_C \lrw \Te_{\dot C} \oplus \Te_{\ddot C} \lrw
\Te_{\dot C \cap \ddot C} \rw 0 .$$Therefore $\X(C)=
\X(\dot C)+\X(\ddot C)-\X(\dot C\cap\ddot C)$ where 
$\X(C)$ is the Euler characteristic, whence $g(C)=g(\dot
C)+g(\ddot C)+
\kappa -1$.
\qed

\begin{rem} \label{union}
{\em Let $C^d_g=\bigcap_{k=1}^{r}
\om(\p^{n_k}, \p^n)$ be a curve in $G(1, n)$, which
defines an incidence scroll in $\p^n$. If we want to
obtain another curve $\dot C$ which is an intersection of
Schubert varieties with the same genus, then $\ddot C$ must
be a line, there is only one generator in common and
$deg(\dot C)=d-1$. 

Let $D^{d_k}_g$ be the directrix curve of
$R^d_g \sub \p^n$ contained in $\p^{n_k}$ for
$k\in \{1,\cdots ,r\}$. Then the directrix curve of $R^{d-1}_g \sub
\p^{n-1}$ in the $\p^{n_k-1}$ has degree $d_k-1$, 
for $k\not=i,j$. This is not so in the other cases
because if $k=i,j$, then the directrix curve in
$\p^{n_k}$ has degree $d_k$ (an easy computation of
Schubert cycles).

We shall refer to this particular degeneration as join
$\p^{n_i}$ and $\p^{n_j}$ and to the inverse as separate
$\p^{n_i}$ and $\p^{n_j}$.}
\end{rem}

\begin{example1}\label{ultima}
{\em The following examples are a good illustration of join
and separate $\p^{n_i}$ and $\p^{n_j}$. 

\begin{enumerate}
\item[$\bullet$]$\be_1 =\{2\, \p^{n-1},
2\, \p^{n}, \p^{2n-3}\} \Longrightarrow R^{4n-6}_{n-2}\sub\p^{2n-1}$.

If we suppose $\p^{n} \cup
\p^{2n-3}=\p^{2n-2}$, i.e.,
$\p^{n} \cap \p^{2n-3}=\p^{n-1}$, then the scroll breaks up into
$R^{2n-2}_{0}\sub\p^{2n-1}$  with base $\{3\, \p^{n-1}, \p^{n}\}$ and
$R^{2n-4}_{0}\sub\p^{2n-3}$  with base $\{3\, \p^{n-2}, \p^{n-1}\}$
(we will see it in section \ref{4}) and
$n-1$ generators in common. It follows that $d=2n-2+2n-4=4n-6$ and $g=n-2$.
\item[$\bullet$]$\be_2 =\{ \p^{n-2},
2\, \p^{n-1}, \p^{n}, \p^{2n-4}\} \Longrightarrow 
R^{4n-9}_{n-3}\sub\p^{2n-2}$; apply \ref{union} to a $\p^{n-1}$
and a 
$\p^{n}$ in $\be_1$ and write $n$ instead of $n-1$.
\item[$\bullet$]$\be_3 =\{ \p^{n-3},
3\, \p^{n-1}, \p^{2n-5}\} \Longrightarrow
R^{4n-12}_{n-4}\sub\p^{2n-3}$; apply \ref{union} to $\p^{n-2}$
and $\p^{n}$ in $\be_2$ and write $n$ instead of $n-1$.
\item[$\bullet$]$\be_4 =\{3\, \p^{n-1},
\p^{n+1}, \p^{2n-3}\} \Longrightarrow
R^{4n-8}_{n-3}\sub\p^{2n-1}$; apply \ref{union} to $2\,
\p^{n-1}$'s in $\be_2$. 
\end{enumerate}
We will use the existence of the previous incidence scrolls
in section \ref{seiss}.
}
\end{example1}

\section{Incidence Rational Scrolls}\label{4}

Let $X=\p(\E)$ be a ruled surface. We say that $X$
is a rational ruled surface if $C\cong \p^1$.
One sees immediately that $X_0=\p^1\times \p^1$ with
its first projection is a rational ruled surface. Each
rational ruled surface is isomorphic to the
ruled surface obtained from $X_0$ by applying a finite
number of elementary transform at
$P_1,\cdots,P_e \in X_0$, with a suitable number $e\geq0$.
We shall write the above expression as $X_e=
elem_{(P_1,\cdots,P_e)}(X_0)$. The integer $e$ is an
invariant of $X_e$. Moreover, for each
$e\geq 0$ there is exactly one (up to isomorphism) rational
ruled surface with invariant $e$, given by $\E
=\Te_{\p^1}\oplus\Te_{\p^1}(-e)$, i.e., $X_e \cong \p(\E)$.

\begin{prop} \label{in}
Let $D\sim aC_o+bf$ be a divisor on $X_e$. Then:
\begin{enumerate}
\item $D$ is very ample $\iff$ $D$ is ample $\iff$ $a>0\, $
and $\, b>ae$;
\item $| D |$ contains an irreducible
nonsingular curve $\iff$
$| D |$ contains an irreducible curve  $\iff$
$a=0, b=1$; or $ a=1,b=0$;
or $a>0, b>ae$; or $e>0,a>0,b=ae$.\qed
\end{enumerate}
\end{prop}

Let $H$ be a very ample divisor on $X_e$ and
let $\Phi_H$ be the closed immersion defined by $| H |$.
Then we know that $\Phi_H(X_e)$ is a scroll if and only 
if $H\equiv C_o+mf$ with $m>e$. To see a more complete 
theory of rational ruled surfaces we refer the reader to 
\cite{hatshor} and \cite{Lanteri}.

The aim of this chapter is to study all the
incidence rational scrolls with base in general
position. For this, we can now formulate a theorem
as follows.

\begin{teo}\label{cla}
Let $H\sim C_o+mf$ be a very ample divisor on $X_e$
and let {\small $\Phi_H:X_e \rw R_0^{2m-e} \sub
\p^{2m-e+1}$} be the closed immersion defined by $| H |$.
Then $R_0^{2m-e}$ is an incidence scroll if and only if
it satisfies one of the following conditions:
\begin{enumerate}
\item $H\sim C_o+(e+1)f$ ($\Phi_H(C_o)=\p^1$);
\item $e=0$ ($\Phi_H(C_o)=C_0^m\sub \p^m$);
\item $e=1$ ($\Phi_H(C_o)=C_0^{m-1}\sub
\p^{m-1}$).
\end{enumerate}
\end{teo}
We will divide its proof into a sequence of propositions.

\subsection*{Rational Scrolls with a Directrix Line}

It is clear that an incidence scroll with 
$\p^1$ as base space is a rational scroll. The following 
proposition gives all rational
normal scrolls with a line as directrix. The affirmation is
not true if normal is deleted from the hypothesis. For
example, $R^3_0 \sub \p^3$, which has a double line, is not an
incidence scroll if the base spaces are in general
position.

\begin{prop}	\label{rec}
In $\p^n$ the incidence scroll with base $\be_n =\{\p^1,
(n-1)\p^{n-2}\}$ is the rational normal scroll of degree
$n-1$ with a line as minimum directrix.
\end{prop}
{\bf Proof.} We proceed by induction in $n$. The proposition
is true in $\p^3$. In this case, the scroll is the
quadric surface in $\p^3$.

Supposing the proposition true in $\p^n$, we prove it for
$n+1$. We can add a hyperplane to the base without
affecting the scroll. Separating $\p^1$ and $\p^{n-1}$, as in
\ref{union}, we obtain a base
$\be_{n+1} =\{\p^1, n\, \p^{n-1}\}$ which defines $R^{n}_0\sub
\p^{n+1}$ as incidence scroll.\qed

\begin{cor} \label{subrec}
For every $e\geq 0$, there is an embedding of $X_e$ as 
incidence rational scroll.
\end{cor}
{\bf Proof.} Use the very ample divisor $H\sim
C_o+(e+1)f$.\qed

\subsection*{Rational Scrolls of General Type}

$X_e$ is said to be of general type if and only if $\E
\cong \Te_{\p^1} \oplus \Te_{\p^1}$ or $\E
\cong \Te_{\p^1} \oplus \Te_{\p^1}(-1)$. Moreover, if $X_e$
is general of degree $d$, then the
minimum degree of the directrix curves is equal to
$\frac{d}{2}$ (respectively,
$\frac{d-1}{2}$) for $d$ even (respectively, odd). On $X_e$,
the family of directrices of degree $m$ has dimension
$d_m=2m-d+1$ (see \cite{Ghione}, pp. 89-90).

\begin{prop}\label{x0x1}
All general rational scrolls are incidence scrolls. 
\end{prop} 
{\bf Proof.}  {\bf e=0:} We will see that in
$\p^{2m+1}$ the incidence scroll generated by $\be
=\{3\,\p^m,\p^{m+1}\}$ is in fact $R^{2m}_0\sub \p^{2m+1}$.
This incidence scroll have directrix curves in each $\p^m$
of degree $d_m$ and a  directrix curve in $\p^{m+1}$ of
degree $d_{m+1}$. Then:
\begin{enumerate}
\item $d_{m+1}$ is the number of lines in $\p^{2m+1}$ which
meet
$4\,\p^{m}$'s. Call this number $z(m)$ i.e.,
$d_{m+1}=z(m)=\Omega(m,2m+1)^4$. 
\item if we take a
hyperplane $\p^{2m}\sub\p^{2m+1}$ through the base
$\p^{m+1}$, then we can calculate $d_m$ as the number of
generators of the scroll which meet $\p^{2m}$ in points
of the base $\p^{m}$. To obtain this number we consider two
cases: either the generators lie in $\p^{2m}$ or they do
not. In the former case, we have the lines in $\p^{2m-3}$
which meet
$4\,
\p^{m-2}$'s therein. The number of these lines is $z(m-2)$.
In the latter case, this number is equal to the number of
lines through $\p^{m}\cap \p^{m+1}$, which meet the other
$2\, \p^{m}$'s. Thus
$d_m=z(m-2)+1$.
\end{enumerate}

Therefore $deg(R)=d_m+d_{m+1}-1=2d_m \Longrightarrow
z(m)=z(m-2)+2 \Longrightarrow d_{m+1}=z(m)=d_m+1$.
But
$z(m)$ is the number of lines meeting
$4\, \p^{m}$'s in $\p^{2m+1}$, i.e., by Schubert calculus,
$z(m)=m+1$. Then
$deg(R)=2m$. Moreover, a scroll of degree $2m$ in $\p^{2m}$
is necessarily rational and  normal.

{\bf e=1:} Since $R_0^{2m-1} \sub \p^{2m}$ is of
general type, it has one minimum directrix curve. Then we can take as
base space a $\p^{m-1}$ which contains it. Choose
$3\,\p^m$'s containing 3 generic directrix curves
$C^{m+1}_0$. This is possible because
$h^0(\Te_{X_1}(C_o+f))=3$. By Proposition \ref{in} and
Bertini's Theorem we know that
$|C_o+f|$ has irreducible curves and that the generic is
irreducible. $\be =\{\p^{m-1},3 \, \p^{m}\}$ is base of
$R_0^{2m-1} \sub
\p^{2m}$. We conclude it from \ref{union}, joining $2\,
\p^m$ of case
$e=0$.\qed

\subsection*{Rational Scrolls which are not Incidence
Scrolls}
\begin{prop}\label{noin}
Let $H\sim C_o+mf$ be a very ample divisor
on $X_e$ with $m>e+1>2$. If the base must be in
general position, then $\Phi_H(X_e)=R_0^{2m-e}\sub
\p^{2m-e+1}$ cannot be obtained by incidences.
\end{prop}
{\bf Proof.} Under the
above hypotheses, we have 
$m+(e+2)(m-e)>4m-2e-1(\iff 1<e<m-1)$.
Then:$$\begin{array}{lr} m+(\eta)(m-e)\leq 4m-2e-1& (*) \\
m+(\eta+1)(m-e)>4m-2e-1& \\
\end{array}$$for a suitable number $\eta $ such that
$1<\eta<e+2$. An argument similar to Proposition \ref{l62}
shows that there are $\eta\,\p^m$'s in the
base. Then $\be=\{ \p^{m-e}, \eta \,\p^{m},
\p^{n_1},
\cdots ,\p^{n_r}\}$ with $m+1\leq n_i\leq 2m-e-1$
is base of an incidence scroll, which is not
rational. In other case, joining base spaces, we
find a rational scroll in $\p^3$, obtained by
Proposition \ref{rec}, which is impossible. If we
apply it successively the inverse degeneration
(separate $\p^{n_i}$ and $\p^{n_j}$) to all the
possible pairs
$(n_i, n_j)$ such that $\p^{n_i}\cap\p^{n_j}=\emptyset$,
then we must obtain our scroll, but we obtain the
rational normal scroll $R^{2m-e}_0\sub
\p^{2m-e+1}$ with a line as minimum directrix curve (i.e.,
$m=e+1$) or the rational scroll of general type (i.e.,
$e\leq 1$). Note that we can project to any 
$\p^{2m-e-k}$, and so join two base spaces,
because the scroll is rational and, in particular,
decomposable, i.e.,
$\p^{n_1} \cap \p^{n_2}=
\emptyset $. If (*) is an equality,
we apply this argument again, with $\be$ replaced by
$\{\p^{m-e},\eta \,
\p^{m}\}$. Then no the rational scroll has a base which
defines it as incidence scroll.
\qed

\begin{cor}\label{c53} 
Let $R=R^{2m-e}_0 \sub \p^{2m-e+1}$ be the projective model
of $X_e$ defined by $H\sim C_o+mf$. For $e\geq 1$, $R$ is an
incidence scroll (with base in general position) if and only
if
$mh^0(\Te_X(C_o))+(m-e)(h^0(\Te_X(C_o+ef)))=2(2m-e+1)-3$.
\qed\end{cor}
\begin{example1}
{\em Consider $R^6_0\sub \p^7$ with a directrix conic, a
three-dimensional family of directrix quartics and a
five-dimensional one of directrix quintics. This
surface is a projective model of $X_2$ by $H \sim C_o+4f$. If
we apply Corollary \ref{c53}, then we see that it is not an
incidence scroll because $4+4(4-2)
\not= 2\cdot 7-3$.}
\end{example1}

From Propositions \ref{rec}, \ref{x0x1} and \ref{noin}, we
obtain the proof of Theorem \ref{cla}. In these
proofs, we have built a base for each incidence rational
scroll in $\p^n$. It is easy notice that if a rational
scroll is an incidence scroll, then there is really only one way
to choose a base, i.e., the sequence $n_1, \cdots, n_r$ is
unique. Table
$1$ contains all the incidence rational scrolls up $\p^{8}$.

\begin{center}
\begin{tabular}{|c||c|c|c|c|c|c||c||c||c|}
\hline
\hline\multicolumn{10}{|c|}{TABLE 1. INCIDENCE RATIONAL
SCROLLS }
\\
\hline
\hline
{ }& \multicolumn{6}{|c||}{$n_i, \,
i=1,\cdots ,6$}&{ }&{ }&{ } \\
\cline{2-7} 
{\scriptsize
Scroll}&{\scriptsize
$1$}&{\scriptsize $2$}&{\scriptsize
$3$}&{\scriptsize $4$}&{\scriptsize
$5$}&{\scriptsize
$6$}&{
\scriptsize Min. Dir.$({\star})$}&{\scriptsize
Normalized}&{\scriptsize
$deg(\d)$}\\
\hline
\hline
{\scriptsize $R^2_0\sub
\p^3$}&{\small
3}&{-}&{-}&{-}&{-}&{-}&{\scriptsize $\p^1\;(\infty^1)$}&
{\tiny $\Te_{\p^1}\oplus\Te_{\p^1}$}&{\small 1}\\
\hline
\hline
{\scriptsize $R^3_0\sub \p^4$}&{\small
1}&{\small 3}&{-}&{-}&{-}&{-}&{\scriptsize
$\p^1\;(1)$}&{\tiny $\Te_{\p^1}\oplus
\Te_{\p^1}(-1)$}&{\small 2}\\
\hline
\hline
{\scriptsize $R^4_0\sub
\p^5$}&{\small
1}&{-}&{\small
4}&{-}&{-}&{-}&{\scriptsize $\p^1\;(1)$}&{\tiny
$\Te_{\p^1}\oplus
\Te_{\p^1}(-2)$}&{\small 3}\\
\hline
{\scriptsize $R^4_0\sub
\p^5$}&{-}&{\small 3}&{\small
1}&{-}&{-}&{-}&{\scriptsize $C^2_0\sub
\p^2\;(\infty^1)$}&{\tiny
$\Te_{\p^1}\oplus\Te_{\p^1}$}&{\small 2}\\
\hline
\hline
{\scriptsize $R^5_0\sub
\p^6$}&{\small
1}&{-}&{-}&{\small
5}&{-}&{-}&{\scriptsize $\p^1\;(1)$}&{\tiny
$\Te_{\p^1}\oplus
\Te_{\p^1}(-3)$}&{\small 4}\\
\hline
{\scriptsize $R^5_0\sub
\p^6$}&{-}&{\small 1}&{\small
3}&{-}&{-}&{-}&{\scriptsize $C^2_0\sub
\p^2\;(1)$}&{\tiny
$\Te_{\p^1}\oplus\Te_{\p^1}(-1)$}&{\small 3}\\
\hline
\hline
{\scriptsize $R^6_0\sub
\p^7$}&{\small
1}&{-}&{-}&{-}&{\small
6}&{-}&{\scriptsize $\p^1\;(1)$}&{\tiny
$\Te_{\p^1}\oplus
\Te_{\p^1}(-4)$}&{\small 5}\\
\hline
{\scriptsize $R^6_0\sub
\p^7$}&{-}&{-}&{\small 3}&{\small
1}&{-}&{-}&{\scriptsize $C^3_0\sub
\p^3\;(\infty^1)$}&{\tiny
$\Te_{\p^1}\oplus\Te_{\p^1}$}&{\small 3}\\
\hline
\hline
{\scriptsize $R^7_0\sub
\p^8$}&{\small
1}&{-}&{-}&{-}&{-}&{\small
7}&{\scriptsize $\p^1\;(1)$}&{\tiny $\Te_{\p^1}\oplus
\Te_{\p^1}(-5)$}&{\small 6}\\
\hline
{\scriptsize $R^7_0\sub
\p^8$}&{-}&{-}&{\small 1}&{\small
3}&{-}&{-}&{\scriptsize $C^3_0\sub
\p^3\;(1)$}&{\tiny
$\Te_{\p^1}\oplus\Te_{\p^1}(-1)$}&{\small 4}\\
\hline
\hline\multicolumn{10}{l}{$^{\star}$ Number of minimum directrix
curves}\\
\end{tabular}
\end{center}

\section{Incidence Elliptic Scrolls}\label{seiss}

Let $X=\p(\E)$ be a ruled surface over a elliptic curve
$C$. Let $\de$ be a divisor on $C$. Then:
\begin{enumerate}
\item $deg(\de) \geq e+2 \iff $ there is a section $D\sim
C_o+\de f$ and $|D|$ has no base points;
\item $deg(\de )\geq e+3 \iff$ $|C_o+\de f|$ is very ample.
\end{enumerate}

If we take $H \sim C_o+\d f$ a very ample
divisor on $X$ with $m=deg(\d )$, then we obtain the closed 
immersion $\Phi_{H} \colon X \lrw R^{2m-e}_1\sub
\p^{2m-e-1}$. We denote the scroll briefly by $R$.

We want to know if such a scroll is an incidence scroll. For
this we need a base $\be$. Since such a base is formed by
linear spaces which contain directrix curves of the scroll,
we begin to study the families of directrix curves in $R$.
The aim of this section is to prove the following
theorem.

\begin{teo}\label{c.e}
Let $H\sim C_o+\d f$ be a very ample divisor on $X$ with
$m=deg(\d )$ and let {\small $\Phi_H:X \rw R_1^{2m-e} \sub
\p^{2m-e-1}$} be the closed immersion defined by $| H |$.
Then $R_1^{2m-e}$ is an incidence scroll (with base in
general position) if and only if
it satisfies one of the following conditions:
\begin{enumerate}
\item $e=-1$ and $m=2$
($\be=\{5 \, \p^2\}$);
\item $\e\sim 0$ and $m=4$ ($\be=\{3\, \p^3, 2\, \p^5\}$);
\item $X$ decomposable, $ 0\leq e\leq 3$ and
$m=e+3$ ($\be=\{\p^2, (e+1)\, \p^{e+2}, (3-e)\,
\p^{e+3}\}$).
\end{enumerate}
\end{teo}

\subsection*{Decomposable incidence elliptic scrolls} 

Let $X=\p(\E)$ be a ruled surface over an elliptic curve $C$,
defined by a normalized bundle $\E \cong\Te_C\oplus
\Te_C(\e)$ and let $H\sim C_o+\d f$ be a very ample divisor
on $X$. Suppose $R$ an incidence scroll with base in
general position.

Consider $e\geq 1$. Since $h^0(\Te_X(C_o))=1,$ there
does not exist another directrix curve of minimum degree. The
following directrix curves of the scroll are linearly
equivalent to $C_o-\e f$ and, by $H$, are directrix
curves of degree $m$. Since
$h^0(\Te_X(C_o-\e f))=e+1
\geq 2$, we know that there are $(e+1)\,
\p^{m-1}$'s in
$\be$, whenever possible ($m-1+(e+1)(m-e-1)\leq 4m-2e-5$).

\begin{lemma1}\label{co65} 
Let $X=\p(\Te_C\oplus \Te_C(\e))$ be a ruled surface over
an elliptic curve $C$ with $e\geq1$and let $H\sim
C_o+\d f$ be a very ample divisor on
$X$. Then
$R$ is an incidence scroll if and only if
$deg(\d)=e+3$ and $1 \leq e \leq 3$. Moreover,
$\be=\{\p^2, (e+1)\, \p^{e+2}, (e-3)\,
\p^{e+3}\}$.
\end{lemma1}
{\bf Proof.} We know that there are a $\p^{m-e-1}$ and
$(e+1)\,\p^{m-1}$'s in $\be$ when $e=1,2,3$. For $e=1$, the
following directrix curves are divisors $C\sim C_o+(P+Q)f,
\,P,Q\in C$. Since $h^0(\Te_X(C_o+(Q-\e)f))=3$, we show
that the generic curve is irreducible. For $m\geq 5$, we
can see that one $\p^{m-2}$, $2\,\p^{m-1}$'s and
one $\p^m$ are in $\be$, i.e., we impose
$4m-8$ independent conditions on $G(1,2m-2)$ so that  
$\be =\{\p^{m-2}, 2\, \p^{m-1}, \p^m,\p^{2m-4}\}$ is base of
an incidence scroll. Such a scroll is of type
$R^{4m-9}_{m-3}\sub \p^{2m-2}$ (see Example \ref{ultima})
which is not elliptic. For $m=4$, $\be =\{\p^2, 2\,
\p^3, 2 \, \p^4\}$ is the base of $R^7_1 \sub \p^6$. 

For $e=2$, we only need another condition, i.e., 
$\be=\{\p^{m-3}, 3\p^{m-1},\p^{2m-5}\}$. Such a scroll is
$R^{4m-12}_{m-4}\sub \p^{2m-3}$ (see Example \ref{ultima}),
which is elliptic if and only if $m=5$.

For $e=3$, we have seen that if $R$ is an
incidence scroll, then a $\p^{m-4}$ and
$4\, \p^{m-1}$'s are in $\be$. If $m=6$, $\be=\{\p^{2}, 4\,
\p^{5}\}$ is base of $R^9_1 \sub \p^8$. For $e\geq
3$ and $m\not= 6$, we have
$m-1+(e+1)(m-e-1)>4m-2e-5$. Consequently, we
cannot take $(e+1)\,
\p^{m-1}$'s as base spaces because these impose
too many conditions on $G(1,2m-e-1)$. Since $m-1<2m-e-3$, we
conclude (as in Proposition \ref{noin}) that $R$ is not an
incidence scroll.\qed

\begin{rem} {\em We have seen that, in $\p^N$, if the base
space of greater dimension is $\p^{N-2}$'s, then it is not
necessary to take the number of linearly independent
directrix curves. From this, we need to show if
$m-1<2m-e-3$.}
\end{rem}

We now evaluate the case $e=0$. There are two choices of the
divisor $\e$. If $\e \not\sim 0$, there are
two choices of normalized $\E$, namely $\Te_C\oplus
\Te_C(\e)$ and $\Te_C\oplus \Te_C(-\e)$. There are exactly 
two choices of $C_o$, both with $C_o^2=0$, namely $C_o$ and
$C_o-\e f$. Since $h^0(\Te_X(C_o))=h^0(\Te_C(C_o-\e f))=1$,
these are the only minimum directrix curves. If $R$ is an
incidence scroll with base in general position, then there
are $2\, \p^{m-1}$'s in $\be$. The following directrix
curves are in linear systems $| C_o+(P+\e) f|$ and $|
C_o+(P-\e) f|, \, P \in C$, where $h^0(\Te_X(C_o+(P+\e)
f))=h^0(\Te_X(C_o+(P-\e) f))=2$. As in the proofs
described above we can see that, if $R$ is an incidence
scroll, there are $2\, \p^m$'s as base spaces. Then $\be
=\{2\, \p^{m-1},2\, \p^m, \p^{2m-3}\}$ and the incidence
scroll has degree $4m-6$ and genus
$m-2$ (see Example \ref{ultima}), which is elliptic
if and only if
$m=3$.

In the other case, $\e\sim 0$, there is a
one-dimensional family of directrix curves of degree
$m$. If $R$ is an incidence scroll, there are $3\,
\p^{m-1}$'s in $\be$. The following
directrix curves are in a linear system
$|C_o+(P+Q)f|$ where
$h^0(\Te_{X}(C_o+(P+Q)f))=4$. By Bertini, we know that it has
irreducible curves and that so is the generic curve. If $R$ is
an incidence scroll, then $m\geq 4$ (no three
$\p^{m-1}$'s suffice and another  $\p^{m-1}$ imposes too
many conditions). For $m\geq 4$, $\be=\{3\,
\p^{m-1},\p^{m+1},\p^{2m-3}\}$ generates a scroll of
genus $m-3$ (see Example
\ref{ultima}), which is elliptic if and only if $m=4$.

\begin{rem}
{\em If we want to extend the study of incidence elliptic
scrolls to incidence scrolls of genus $g$, then we must
consider the following properties.
\begin{enumerate}
\item If the scroll has a unique minimum directrix curve,
then the space which contains it is in
$\be$.
\item If the scroll has two minimum directrix curves, then the
correspondent spaces are in
$\be$.
\item If the scroll has an one-dimensional family of minimum
directrix curves, then the base has
three spaces which contain three of these curves.
\item If $m-g+(e+2-g)(m-e-g)>4m-2e-4g-1$, then $R$ is
not an incidence scroll ($g=0,1$). 
\end{enumerate}} 
\end{rem}

\subsection*{Indecomposable incidence elliptic scrolls} 
Let $X=\p(\E)$ be a ruled surface
over an elliptic curve $C$, corresponding to an
indecomposable $\E$. Then $e=0,-1$ and there
is exactly one ruled surface over $C$ for each of these
values of $e$. Let $H\sim C_o+\d f$ be a divisor with 
$deg(\d)=m\geq e+3$, which defines 
$\Phi_{H}\colon X\lrw R^{2m-e}_1\sub \p^{2m-e-1}$.
We want to know if such a scroll is an incidence
scroll. For this, we need to find a suitable
number of base spaces in $\p^{2m-e-1}$ whose 
dimension is $\leq 2m-e-3$ 

If $e=0$ then $\E$ is an extension of $\Te_C$ by $\Te_C$
i.e. given by an exact sequence $0\rw \Te_C \lrw \E \lrw
\Te_C \rw 0$. Since $C_0^2=0$ and $X$ is
indecomposable, there is a unique directrix of minimum
degree equal to $m$. For any $P\in C$,
$h^0(\Te_X(C_o+Pf))=2$. Then we have an one-dimensional
family $|C_o+Pf|$ of directrix curves of degree
$m+1$, for any $P\in C$. Suppose that $R$ is an
incidence scroll with base in general position. An
easy computation shows that $\be=\{\p^{m-1}, 2\,
\p^m, \p^{n_1},
\cdots, \p^{n_r}\}$ where $m\leq n_i \leq 2m-3$ and
$\sum_{i=1}^r(2m-2-n_i)=m$. Then we can obtain
another incidence elliptic scroll in
$\p^{2m}$  with $e=1$ and base $\dot
\be=\{\p^{m-1},\p^m, \p^{m+1}, \p^{n_1+1},\cdots,
\p^{n_r+1}\}$ (separate
$\p^{m-1}$ and $\p^m$ in $\be$). This happens when 
$m=3, \, n_1=2$ and $n_2=3$, which contradicts the
hypothesis $m\leq n_i$. Hence there are no indecomposable
elliptic scrolls with
$e=0$ defined by incidences.

If $e=-1$ then we have an exact sequence $0\rw \Te_C \lrw
\E\lrw\Te_C(P) \rw 0$. Then $\Phi_{H}(C_o)$ is a directrix
curve of degree $m+1$. If $D$ is a section with $D^2=1$,
then $D\sim C_o+(\delta -\e)f$ with $deg(\delta)=1$. In
fact, the sections $C_o$ with $C_o^2=1$ form a
one-dimensional algebraic family parametrized by $C$ such
that no two of them are linearly equivalent. If $R$ is an
incidence scroll, then $\be=\{2\, \p^m, \p^{n_1},
\cdots, \p^{n_r}\}$ with $m\leq n_i \leq 2m-2$ and
$\sum_{i=1}^r(2m-n_i-1)=2m-1$. We can obtain
another incidence elliptic scroll in $\p^{2m+1}$ 
with $\e \not\sim 0$ and $e=0$ (separate $2\,
\p^m$'s in
$\be$). Consequently, there is only one incidence scroll
with $e=-1$. This happens when 
$m=n_i=2$ and $r=3$. So $\Phi_H(X)=R^5_1\sub
\p^4$ is an incidence scroll with base $\be =\{5\, \p^2\}$.

In Table $2$ we have compiled all incidence elliptic
scrolls. The interest of this table is that each scroll
is a degenerate scroll of the remainder of them.

\begin{center}
\begin{tabular}{|c||c|c|c|c||c||c||c|}
\hline
\hline\multicolumn{8}{|c|}{ TABLE 2. INCIDENCE
ELLIPTIC SCROLLS}\\
\hline
\hline
{ }& \multicolumn{4}{|c||}{\scriptsize $n_i, \,
i=1,\cdots ,4$}&{ }&{ }&{ } \\
\cline{2-5} 
{\small Scroll}&{\small
$2$}&{\small $3$}&{\small $4$}&{\small $5$}&
{\small Min. Dir.}&{\small Normalized}&{\small
$deg(\d)$}\\
\hline
\hline
{\small $R^5_1\sub \p^4$}&{\small 5}&{\small
-}&{\small -}&{\small -}&{\small
$C^3_1\sub
\p^2$}&{\tiny $\E \in
Ext^1(\Te_{C}(P),\Te_{C})$}&{\small 2}\\
\hline
\hline
{\small $R^6_1\sub \p^5$}&{\small 2}&{\small
3}&{\small -}&{\small -}&{\small
$C^3_1\sub
\p^2$}&{\tiny $\Te_{C}\oplus\Te_{C}(\e)(\e\not\sim
0$)}&{\small 3}\\
\hline
\hline
{\small $R^7_1\sub \p^6$}&{\small 1}&{\small
2}&{\small 2}&{\small -}&{\small
$C^3_1\sub
\p^2$}&{\tiny $\Te_{C}\oplus\Te_{C}(-P)$}&{\small 4}\\
\hline
\hline
{\small $R^8_1\sub \p^7$}&{\small 1}&{\small
-}&{\small 3}&{\small 1}&{\small
$C^3_1\sub
\p^2$}&{\tiny $\Te_{C}\oplus\Te_{C}(-P-Q)$}&{\small 5}\\
\hline
{\small $R^8_1\sub
\p^7$}&{\small -}&{\small 3}&{\small -}&{\small 2}&{\small
$C^4_1\sub
\p^3$}& {\tiny $\Te_{C}\oplus\Te_{C}$}&{\small 4}\\
\hline
\hline
{\small $R^9_1\sub \p^8$}&{\small 1}&{\small
-}&{\small -}&{\small 4}&{\small
$C^3_1\sub
\p^2$}&{\tiny $\Te_{C}\oplus\Te_{C}(-P-Q-R)$}&{\small 6}\\
\hline
\hline
\end{tabular}
\end{center}

\begin{flushleft}
{\bf Acknowledgements:}
Rosa Cid-Mu$\tilde n$oz was supported by a 
fellowship of Xunta de Galicia (Spain).
\end{flushleft}

\end{document}